\documentclass[10pt]{article}
\usepackage{graphicx}
\usepackage{amssymb}
\newtheorem{theorem}{Theorem}

\newtheorem{lemma}[theorem]{Lemma}
\newtheorem{corollary}[theorem]{Corollary}

\input tcilatex

\begin{document}

\begin{titlepage}

\vskip0.5truecm

\vskip1.0truecm

\begin{center}

{\LARGE \bf Uniform bounds for diffeomorphisms of the torus and a conjecture of P. Boyland}

\end{center}

\vskip  0.4truecm

\centerline {{\large Salvador Addas-Zanata}}

\vskip 0.2truecm

\centerline { {\sl Instituto de Matem\'atica e Estat\'\i stica }}
\centerline {{\sl Universidade de S\~ao Paulo}}
\centerline {{\sl Rua do Mat\~ao 1010, Cidade Universit\'aria,}} 
\centerline {{\sl 05508-090 S\~ao Paulo, SP, Brazil}}
 
\vskip 0.7truecm

\begin{abstract}

We consider $C^{1+\epsilon}$ diffeomorphisms of the 
torus, denoted $f,$ homotopic to the identity and whose rotation sets have interior. We give some 
uniform bounds on the displacement of points in the plane under
iterates of a lift of $f,$ relative to vectors in the 
boundary of the rotation set and we use these estimates in order to prove that if 
such a diffeomorphism $f$ preserves area, then the rotation vector of 
the area measure is  an interior point of the rotation set. This settles a 
strong version of a conjecture proposed by P. Boyland. We also present some new 
results on the realization of extremal points of the rotation set by compact
$f$-invariant subsets of the torus.

\end{abstract} 

\vskip 0.3truecm

\vskip 2.0truecm


\vskip 0.8truecm

\noindent{\bf e-mail:} sazanata@ime.usp.br

\vskip 1.0truecm


\vfill
\hrule
\noindent{\footnotesize{The author is partially supported 
by CNPq, grant: 303127/2012-0}}

\end{titlepage}

\baselineskip=6.2mm

\section{Introduction and main results}

The main motivation for this paper is to study how rigid is the displacement
of points in the plane under the action of a lift of a homeomorphism of the
two dimensional torus homotopic to the identity (more precise explanations
will be given below). The similar problem for an orientation preserving
homeomorphism of the circle was already studied by H. Poincar\'e. He proved
that given an orientation preserving circle homeomorphism $f:S^1\rightarrow
S^1$ and a lift of $f$ to the real line, denoted $\widetilde{f}:{\rm I%
\negthinspace R\rightarrow I\negthinspace R,}$ there exists a number $\omega
\in {\rm I\negthinspace R,}$ called the rotation number of $\widetilde{f},$
such that 
$$
\left| \widetilde{f}^n(\widetilde{x})-\widetilde{x}-n.\omega \right| <2, 
\text{ for all }\widetilde{x}\in {\rm I\negthinspace R}\text{ and any
integer }n>0. 
$$

The situation for homeomorphisms of the torus is more complicated. In
general there is no such $\omega $ as above and some points may not even
have a rotation vector, the generalization of rotation number to this new
setting. In order to make things precise and to present our main results and
some motivation, a few definitions are necessary:

{\bf Basic notation and some definitions:}

\begin{enumerate}
\item  Let ${\rm T^2}={\rm I}\negthinspace {\rm R^2}/{\rm Z\negthinspace
\negthinspace Z^2}$ be the flat torus and let $p:{\rm I}\negthinspace {\rm %
R^2}\longrightarrow {\rm T^2}$ be the associated covering map. Coordinates
are denoted as $(\widetilde{x},\widetilde{y})\in {\rm I}\negthinspace {\rm %
R^2}$ and $(x,y)\in {\rm T^2.}$

\item  Let $Diff_0^{1+\epsilon }({\rm T^2})$ be the set of $C^{1+\epsilon }$
(for some $\epsilon >0$) diffeomorphisms of the torus homotopic to the
identity and let $Diff_0^{1+\epsilon }({\rm I\negthinspace R^2})$ be the set
of lifts of elements from $Diff_0^{1+\epsilon }({\rm T^2})$ to the plane.
Maps from $Diff_0^{1+\epsilon }({\rm T^2})$ are denoted $f$ and their lifts
to the plane are denoted $\widetilde{f}.$ By $Diff_0^0({\rm T^2})$ we mean
the set of homeomorphisms of ${\rm T^2}$ homotopic to the identity and $%
Diff_0^0({\rm I\negthinspace R^2})$ is the set of lifts of elements from $%
Diff_0^0({\rm T^2})$ to the plane. In this $C^0$-setting, maps of the torus
are also denoted $f$ and their lifts to the plane are denoted $\widetilde{f}.
$

\item  Let $p_{1,2}:{\rm I}\negthinspace {\rm R^2}\longrightarrow {\rm I}%
\negthinspace {\rm R}$ be the standard projections; $p_1(\tilde x,\tilde
y)=\tilde x$ and $p_2(\tilde x,\tilde y)=\tilde y$. 

\item  Given $f\in Diff_0^0({\rm T^2})$ and a lift $\widetilde{f}\in
Diff_0^0({\rm I}\negthinspace {\rm R^2}),$ the so called rotation set of $
\widetilde{f},$ $\rho (\widetilde{f}),$ can be defined as follows (see \cite
{misiu}): 
\begin{equation}
\label{rotsetident}\rho (\widetilde{f})=\bigcap_{{\ 
\begin{array}{c}
i\geq 1 \\ 
\end{array}
}}\overline{\bigcup_{{\ 
\begin{array}{c}
n\geq i \\ 
\end{array}
}}\left\{ \frac{\widetilde{f}^n(\widetilde{z})-\widetilde{z}}n:\widetilde{z}%
\in {\rm I}\negthinspace {\rm R^2}\right\} }
\end{equation}

This set is a compact convex subset of ${\rm I\negthinspace R^2}$ (see \cite
{misiu}), and it was proved in \cite{franksrat} and \cite{misiu} that all
points in its interior are realized by compact $f$-invariant subsets of $%
{\rm T^2,}$ which can be chosen as periodic orbits in the rational case. By
saying that some vector $\rho \in \rho (\widetilde{f})$ is realized by a
compact $f$-invariant set, we mean that there exists a compact $f$-invariant
subset $K\subset {\rm T^2}$ such that for all $z\in K$ and any $\widetilde{z}%
\in p^{-1}(z)$ 
\begin{equation}
\label{deffrotvect}\stackunder{n\rightarrow \infty }{\lim }\frac{\widetilde{f%
}^n(\widetilde{z})-\widetilde{z}}n=\rho .
\end{equation}
Moreover, the above limit, whenever it exists, is called the rotation vector
of the point $z,$ denoted $\rho (z)$$.$
\end{enumerate}

As the rotation set is a compact convex subset of the plane, there are three
possibilities for its shape:

\begin{enumerate}
\item  it is a point;

\item  it is a linear segment;

\item  it has interior;
\end{enumerate}

An important problem in this set up is to decide which subsets can be
realized as rotation sets of homeomorphisms of the torus homotopic to the
identity. For instance, with a simple rotation, all points can be realized.
Some linear segments can be realized, for others it is not know. And what
about the case when the rotation set has interior. Which sets can be
realized? Rational polygons \cite{kwa1} can, but what else? We do not
consider this problem, but we refer to \cite{frankmisiu} and \cite{kwa2}.%

In the first possibility above, F\'abio Tal and Andr\'es Koropecki \cite
{talkoro} presented an example of an area preserving $C^\infty $
diffeomorphism of the torus homotopic to the identity, denoted $f,$ which
has a lift $\widetilde{f}$ to the plane such that $\rho (\widetilde{f}%
)=\{0\} $ and some points in the plane have unbounded orbits in every
direction. In particular, there exists a point $\widetilde{x}_0\in {\rm I%
\negthinspace R^2} $ such that 
$$
\left| \widetilde{f}^n(\widetilde{x}_0)-\widetilde{x}_0-n.0\right| \text{ is
unbounded with }n>0. 
$$
This type of behavior is usually called sub-linear displacement because,
although there are unbounded $\widetilde{f}$-orbits in the plane, this
behavior is not captured by the rotation set.

Related to the second possibility for the shape of the rotation set, Pablo
Davalos \cite{davalos} analyzed the following situation: Assume $f:{\rm %
T^2\rightarrow T^2}$ is a homeomorphism of the torus homotopic to the
identity and $\widetilde{f}:{\rm I\negthinspace R^2\rightarrow I%
\negthinspace R^2}$ is a lift of $f$ such that some linear segment $
\overline{AB}$ is contained in the boundary of $\rho (\widetilde{f})$ for
some $A,B$ rational vectors. He considered two situations:

\begin{itemize}
\item  $\rho (\widetilde{f})=\overline{AB};$

\item  $\rho (\widetilde{f})$ has interior;
\end{itemize}

In the first case, let $\overrightarrow{v}^{\perp }$ be a unit vector
orthogonal to $\overline{AB}$ with any of the two possible orientations and
in the second, let $\overrightarrow{v}^{\perp }$ be the unit vector
orthogonal to $\overline{AB}$ such that $-\overrightarrow{v}^{\perp }$points
towards $\rho (\widetilde{f}).$ Then Davalos proved the following:

\begin{description}
\item[{Theorem [Davalos]}]  : There exists a number $M>0$ such that 
$$
\left\langle \widetilde{f}^n(\widetilde{x})-\widetilde{x}-n.A,
\overrightarrow{v}^{\perp }\right\rangle \leq M,\text{ for all }\widetilde{x}%
\in {\rm I\negthinspace R^2}\text{ and any integer }n>0. 
$$
\end{description}

Our main result is similar to the above one, but it deals with all the
possible situations when $\rho (\widetilde{f})$ has interior. As our methods
rely on some results from \cite{c1epsilon}, we need a stronger hypothesis,
namely we assume $f\in Diff_0^{1+\epsilon }({\rm T^2}).$

In order to state our main results, let us introduce a little more notation:
Given a compact convex subset $K\subset {\rm I\negthinspace R^2,}$ for every 
$\alpha \in \partial K,$ there exists a straight line $r$ containing $\alpha 
$ such that $K\subset r\cup \{$one connected component of $r^c\}.$ This line
is called a supporting line at $\alpha .$ For instance, in case $\alpha $ is
a vertex, there are infinitely many supporting lines at $\alpha .$

\begin{theorem}
\label{main1}: Let $f\in Diff_0^{1+\epsilon }({\rm T^2})$ be such that $\rho
(\widetilde{f})$ has interior. Then, there exists a number $M_f>0$ such that
for any $\omega \in \partial \rho (\widetilde{f})$ and any supporting line $r
$ at $\omega ,$ if $\overrightarrow{v}^{\perp }$ is the unitary vector
orthogonal to $r,$ pointing towards the connected component of $r^c$ which
does not intersect $\rho (\widetilde{f}),$ then%
$$
\left\langle \widetilde{f}^n(\widetilde{x})-\widetilde{x}-n.\omega ,
\overrightarrow{v}^{\perp }\right\rangle \leq M_f,\text{ for all }\widetilde{%
x}\in {\rm I\negthinspace R^2}\text{ and any integer }n>0. 
$$
\end{theorem}

{\bf Remarks:}

\begin{itemize}
\item  Our proof will show that $M_f$ can be precisely computed from $f$ and
moreover, the same number works for any map in $Diff_0^{1+\epsilon }({\rm T^2%
})$ sufficiently $C^1$-close to $f;$

\item  This theorem may be used as a tool to numerically estimate rotation
sets. For instance if one is considering a family of maps $f_t\in
Diff_0^{1+\epsilon }({\rm T^2})$ an interesting problem connected to our
result is to study how and when $\rho (\widetilde{f}_t)$ changes as $t$
varies;
\end{itemize}

As a corollary of the above result, we prove a stronger version of Boyland's
conjecture in the torus case:

\begin{theorem}
\label{main2}: Let $f\in Diff_0^{1+\epsilon }({\rm T^2})$ be a Lebesgue
measure preserving diffeomorphism such that $\rho (\widetilde{f})$ has
interior. Then the rotation vector of the Lebesgue measure is an interior
point of $\rho (\widetilde{f}).$
\end{theorem}

Remember that the rotation vector of the Lebesgue measure is defined as: 
\begin{equation}
\label{rotvector}\rho (Leb)\stackrel{def.}{=}\int_{{\rm T^2}}\phi (x)dLeb, 
\end{equation}
where $\phi :{\rm T^2\rightarrow I\negthinspace R^2}$ is the displacement
function given by $\phi (x)=\widetilde{f}(\widetilde{x})-\widetilde{x},$ for
any $\widetilde{x}\in p^{-1}(x).$ In general, if we denote by

$$
M_{inv}(f)=\{\text{subset of all }f\text{-invariant Borel probability
measures in }{\rm T^2}\}, 
$$
then for any $\mu \in M_{inv}(f),$ we define the rotation vector of $\mu ,$ $%
\rho (\mu ),$ as 
$$
\rho (\mu )\stackrel{def.}{=}\int_{{\rm T^2}}\phi (x)d\mu . 
$$
These definitions are clearly motivated by Birkhoff's ergodic theorem, since
for every $x\in {\rm T^2}$ and any integer $n>0,$%
$$
\frac 1n\stackrel{n-1}{\stackunder{i=1}{\sum }}\phi \circ f^i(x)=\frac{ 
\widetilde{f}^n(\widetilde{x})-\widetilde{x}}n,\text{ for any }\widetilde{x}%
\in p^{-1}(x). 
$$
So, given $\mu \in M_{inv}(f),$ for $\mu $ $a.e.$ $x\in {\rm T^2,}$
Birkhoff's ergodic theorem implies that the following limit exists%
$$
\stackunder{n\rightarrow \infty }{\lim }\frac 1n\stackrel{n-1}{\stackunder{%
i=1}{\sum }}\phi \circ f^i(x)=\stackunder{n\rightarrow \infty }{\lim }\frac{ 
\widetilde{f}^n(\widetilde{x})-\widetilde{x}}n=\rho (x)\text{ (the rotation
vector of }x\text{)} 
$$
and 
$$
\stackunder{\rm T^2}{\int }\rho (x)d\mu =\int_{{\rm T^2}}\phi (x)d\mu =\rho
(\mu ). 
$$

One last remark about theorem \ref{main2} is the following: the original
problem posed by P. Boyland was to prove that if $interior(\rho (\widetilde{f%
}))\neq \emptyset $ and $\rho (Leb)=(0,0),$ then $(0,0)\in interior(\rho ( 
\widetilde{f})),$ but in the homeomorphism setting. A proof of this result
in this $C^0$-setting was obtained by F\'abio Tal \cite{fabioboy}.

The next result is another easy corollary of theorem \ref{main1} and lemma 
\ref{fabio}. Before stating it, we have to define a few more concepts. Let $%
K\subset {\rm I\negthinspace R^2}$ be a compact and convex subset. We say
that some point $z\in K$ is an extremal point if, whenever $z$ is the convex
combination of two other points $z_1,z_2\in K,$ then either $z=z_1$ or $%
z=z_2.$ Clearly, extremal points are always in the boundary of $K.$ We say
that $z\in K$ is a vertex if $z$ is an extremal point and there are at least
two, which implies infinitely many, supporting lines at $z.$

\begin{corollary}
\label{boundary}: Let $f\in Diff_0^{1+\epsilon }({\rm T^2})$ be such that $%
\rho (\widetilde{f})$ has interior. Suppose for some $\mu \in M_{inv}(f),$ $%
\rho (\mu )\in \partial \rho (\widetilde{f})$ is a vertex. Then, $supp(\mu )$
is a compact $f$-invariant set which realizes the rotation vector $\rho (\mu
).$ Moreover, there exists $M_\mu >0$ such that for every $x\in supp(\mu ),$
for any $\widetilde{x}\in p^{-1}(x)$ and any integer $n>0,$ $\left\| 
\widetilde{f}^n(\widetilde{x})-\widetilde{x}-n.\rho (\mu )\right\| <M_\mu ,$
that is, there is no sub-linear displacement in $supp(\mu ).$ In case $\rho
(\mu )$ is an extremal point, but not a vertex, if we assume that the
intersection of the (unique) supporting line at $\rho (\mu )$ with $\rho (
\widetilde{f})$ is just $\rho (\mu ),$ then $supp(\mu )$ also realizes the
rotation vector $\rho (\mu ).$ But in this case there may be sub-linear
displacement in $supp(\mu ).$
\end{corollary}

{\bf Remarks:}

\begin{itemize}
\item  in the general case when $\rho (\mu )$ is an extremal point, we do
not know if the above corollary holds;

\item  It was proved by Franks \cite{bordoreal} that rational extremal
points of the rotation set are realized by periodic orbits, but for general
extremal points, the problem was open.
\end{itemize}

In an ongoing work with Andre de Carvalho we are generalizing some results
from \cite{c1epsilon} to other surfaces. After that, using the methods from
this paper, we plan to prove a version of theorem \ref{main2} to surfaces of 
$genus\geq 2.$

For homeomorphisms of the torus homotopic to Dehn twists, results analog to
theorem \ref{main1} and \ref{main2} were proved in \cite{eufabra}.

This paper is organized as follows. In the second section we present a
result from \cite{c1epsilon} important for us and an idea of the proof of
theorem 1 in an easy case. In the third section we prove some auxiliary
lemmas and after that, we prove our main theorems.

\section{An important result and some ideas on the proofs}

In \cite{c1epsilon}, we considered diffeomorphisms $f\in Diff_0^{1+\epsilon
}({\rm T^2})$ which preserve area. But the preservation of area is not
necessary to prove the following result, whose proof is contained in the
proof of theorem 6 of \cite{c1epsilon}.

\begin{theorem}
\label{casogeral}: Suppose $f$ belongs to $Diff_0^{1+\epsilon }({\rm T^2})$
and $(0,0)\in int(\rho (\widetilde{f})).$ Then, $f$ has a hyperbolic
periodic saddle point $Q\in {\rm T^2}$ such that any $\widetilde{Q}\in
p^{-1}(Q)$ is $\widetilde{f}$-periodic and for any pair of integers $(a,b),$ 
$W^u(\widetilde{Q})\pitchfork W^s(\widetilde{Q}+(a,b)).$
\end{theorem}

\vskip 0.2truecm

{\bf Remarks:}

\begin{enumerate}
\item  Clearly, the rotation vector of $Q$ is $(0,0).$

\item  By saying that $W^u(\widetilde{Q})\pitchfork W^s(\widetilde{Q}+(a,b))$
we mean that they have a topologically transverse intersection, which of
course is not necessarily $C^1$ transversal. See figure 1 for a picture
which clarifies this. For a precise explanation, see definition 9 (right
before the statement of lemma 1) of \cite{c1epsilon}.

\item  In the proof of theorem \ref{casogeral}, we obtain a $C^1$-transverse
intersection at least when $(a,b)=(0,0).$
\end{enumerate}

The converse of this result is also true, namely if some map $f\in
Diff_0^{1+\epsilon }({\rm T^2})$ has a hyperbolic periodic point $\widetilde{%
Q}$ such that for three non collinear integer vectors $%
(a_1,b_1),(a_2,b_2),(a_3,b_3)$ we have $W^u(\widetilde{Q})\pitchfork W^s( 
\widetilde{Q}+(a_i,b_i)),$ for $i=1,2,3$ and $(0,0)$ belongs to the convex
hull of $\{(a_1,b_1),(a_2,b_2),(a_3,b_3)\},$ then $(0,0)\in int(\rho ( 
\widetilde{f})).$ This follows from the following: The fact that $W^u( 
\widetilde{Q})\pitchfork W^s(\widetilde{Q}+(a_i,b_i))$ implies that we can
produce a topological horseshoe at $Q\in {\rm T^2}$ such that for some
sequence in the symbolic dynamics (one corresponding to points visiting only
one particular rectangle in the horseshoe), there is a periodic orbit for $f$
whose rotation vector is $\left( \frac{a_i}{N_i},\frac{b_i}{N_i}\right) ,$
for some integer $N_i>0.$ And so, 
$$
(0,0)\in interior\text{ }of\text{ }the\text{ }Conv.Hull\{\left( \frac{a_1}{%
N_1},\frac{b_1}{N_1}\right) ,\left( \frac{a_2}{N_2},\frac{b_2}{N_2}\right)
,\left( \frac{a_3}{N_3},\frac{b_3}{N_3}\right) \}, 
$$
which is contained in the interior of $\rho (\widetilde{f})$ because of its
convexity.

The argument used to prove theorem \ref{main1} can be summarized as follows
in the specific situation when $(0,0)\in int(\rho (\widetilde{f})),$ $\omega
=$$(0,1)\in \partial \rho (\widetilde{f})$ and there is a horizontal
supporting line denoted $r$ at $(0,1).$ This is clearly not a general
setting: Both the point $\omega \in \partial \rho (\widetilde{f})$ and the
direction of the supporting line may be irrational, but it is illustrative
of the general strategy.

Let $\widetilde{Q}\in {\rm I}\negthinspace 
{\rm R^2}$ be a hyperbolic periodic point for $\widetilde{f}$ as in theorem 
\ref{casogeral} which by remark 3 after it, has a $C^1$-transverse
homoclinic intersection. Without loss of generality, we can assume that $
\widetilde{Q}$ is fixed, otherwise we consider the map $\widetilde{g}= 
\widetilde{f}^{n_{\widetilde{Q}}},$ where $n_{\widetilde{Q}}$ is the period
of $\widetilde{Q}$ (maybe twice the period if the eigenvalues at $Q$ are
negative). The rotation set changes as $\rho (\widetilde{g})=n_{\widetilde{Q}%
}.\rho (\widetilde{f}).$ So $(0,0)\in int(\rho (\widetilde{g})),$ $(0,n_{ 
\widetilde{Q}})\in \partial \rho (\widetilde{g})$ and there is a horizontal
supporting line denoted $r^{\prime }$ at $(0,n_{\widetilde{Q}}).$ In the
beginning of the proof of theorem \ref{main1} we show that the statement of
the theorem holds for $\widetilde{f},$ if and only if, it holds for $
\widetilde{g},$ which is actually something very easy to prove. So, let us
assume that $n_{\widetilde{Q}}=1.$

The existence of such a point $\widetilde{Q}$ as above implies that there
are arbitrarily small topological rectangles $D_{\widetilde{Q}}\subset {\rm I%
}\negthinspace 
{\rm R^2}$ such that:

\begin{equation}
\label{defreten2l} 
\begin{array}{c}
\widetilde{Q}\text{ is a vertex of }D_{\widetilde{Q}}\text{ and the sides of 
}D_{\widetilde{Q}}, \text{denoted }\alpha _{\widetilde{Q}},\beta _{ 
\widetilde{Q}},\gamma _{\widetilde{Q}}\text{ and }\delta _{\widetilde{Q}} \\ 
\text{ are contained in }W^s(\widetilde{Q}),W^u(\widetilde{Q}),W^s( 
\widetilde{Q})\ \text{and}\ W^u(\widetilde{Q})\ \text{respectively,} 
\end{array}
\end{equation}

see figure 2. As $D_{\widetilde{Q}}$ is arbitrarily small, we can assume
that 
$$
D_{\widetilde{Q}}\cap (D_{\widetilde{Q}}+(a,b))=\emptyset ,\text{ for all
integer pairs }(a,b)\neq (0,0), 
$$
which means that $p(D_{\widetilde{Q}})\subset {\rm T^2}$ is also a
topological rectangle.

Moreover, there exists an integer 
\begin{equation}
\label{defN2lin}N^{\prime \prime }>0\text{ such that for all }n\geq
N^{\prime \prime }\ \text{we have:} 
\end{equation}

\begin{enumerate}
\item  $\widetilde{f}^n(\beta _{\widetilde{Q}})$ and $\widetilde{f}^n(\delta
_{\widetilde{Q}})$ have topologically transverse intersections with $\alpha
_{\widetilde{Q}}+(0,1),\gamma _{\widetilde{Q}}+(0,1)$ and with $\alpha _{
\widetilde{Q}}+(1,0),\gamma _{\widetilde{Q}}+(1,0);$

\item  $\widetilde{f}^n(\gamma _{\widetilde{Q}})\subset \alpha _{\widetilde{Q%
}};$
\end{enumerate}

Now we construct a closed path connected set $\theta \subset {\rm I}%
\negthinspace 
{\rm R^2}$ such that: %
%
%
%

\begin{enumerate}
\item  $\theta =\theta +(1,0);$

\item  $\theta $ contains $D_{\widetilde{Q}}+i(1,0),$ for all integers $i;$

\item  $\theta $ contains two compact simple arcs $\eta _1$ and $\eta _2$ of
the following form: The arc $\eta _1$ starts at $\widetilde{Q},$ goes
through $\widetilde{f}^{N^{\prime \prime }}(\beta _{\widetilde{Q}})$ until
it crosses $\alpha _{\widetilde{Q}}+(1,0)$ and $\gamma _{\widetilde{Q}%
}+(1,0).$ The arc $\eta _2$ starts at $\widetilde{f}^{N^{\prime \prime
}}(\delta _{\widetilde{Q}}\cap \alpha _{\widetilde{Q}}),$ goes through $
\widetilde{f}^{N^{\prime \prime }}(\delta _{\widetilde{Q}})$ until it
crosses $\alpha _{\widetilde{Q}}+(1,0)$ and $\gamma _{\widetilde{Q}}+(1,0),$
see figure 3;

\item  clearly, $\theta $ contains $\eta _{1(or\text{ }2)}+i(1,0),$ for all
integers $i;$

\item  $\theta $ is bounded in the $(0,1)$ direction, that is, $\theta $ is
contained between two straight lines, both parallel to $(1,0),$ and the
distance between them is denoted $d_{(1,0)};$
\end{enumerate}

Now, assume that the uniform bound in the statement of theorem \ref{main1}
does not hold. This means that for every $M>0,$ there exists $\widetilde{x}%
_M\in {\rm I}\negthinspace {\rm R^2}$ and an integer $n_M>0,$ $n_M\stackrel{%
M\rightarrow \infty }{\rightarrow }\infty ,$ such that 
$$
p_2\circ \widetilde{f}^{n_M}(\widetilde{x}_M)-p_2(\widetilde{x}_M)-n_M>M. 
$$

If we choose a sufficiently large $M>0$ and the point $\widetilde{x}_M$
below $\theta $ satisfying $dist.(\widetilde{x}_M,\theta )\leq
2+2.d_{(1,0)}, $ then we get that 
$$
\widetilde{f}^{n_M}(\theta )\text{ intersects }\theta +(0,n_M+\left\lfloor
M-4-4.d_{(1,0)}\right\rfloor ). 
$$
More precisely, for some integer $a,$%
$$
\widetilde{f}^{n_M+N^{\prime \prime }}(D_{\widetilde{Q}})\cap \left( D_{ 
\widetilde{Q}}+(a,n_M+\left\lfloor M-5-4.d_{(1,0)}\right\rfloor )\right) 
$$
contains a connected topological rectangle $\widetilde{R}^{*}$ as in figure
4. So, there is a topological horseshoe in $D_Q=p(D_{\widetilde{Q}})\subset 
{\rm T^2}$ and in particular, this topological horseshoe has a point which
is fixed under iterates of $f^{n_M+N^{\prime \prime }}$ and this point
belongs to $p(\widetilde{R}^{*})\subset D_Q$. So, it has a rotation vector
whose second coordinate is equal to 
$$
\frac{n_M+\left\lfloor M-5-4.d_{(1,0)}\right\rfloor }{n_M+N^{\prime \prime }}%
, 
$$
which is larger than one, if $M>0$ is sufficiently large. So we produced a
point whose rotation vector belongs to the connected component of $r^c$
which does not intersect the rotation set. This contradiction proves the
theorem.

\section{Proofs}

In the first subsection, we prove some auxiliary results.

\subsection{Auxiliary results}

In the next lemma we are going to produce, for every possible direction $
\overrightarrow{v}$$,$ an unbounded closed connected set $\theta _{
\overrightarrow{v}}\subset {\rm I}\negthinspace {\rm R^2}$ which separates
the plane into two special unbounded connected components (maybe there are
other components in the complement of $\theta _{\overrightarrow{v}}),$ by
concatenating integer translates of appropriate pieces of the stable and
unstable manifolds of the hyperbolic $\widetilde{f}$-periodic point $
\widetilde{Q}$ given in theorem \ref{casogeral} (in the applications, the
direction $\overrightarrow{v}$ is that of the supporting line at the
rotation vector $\omega $ in the boundary of $\rho (\widetilde{f})$ we are
considering). The subset $\theta _{\overrightarrow{v}}$ is a general version
of the set $\theta $ considered in the previous section.

As we already explained, for any $\widetilde{Q}\in p^{-1}(Q),$ where $Q$ is
given in theorem \ref{casogeral}, there are arbitrarily small topological
rectangles $D_{\widetilde{Q}}\subset {\rm I}\negthinspace 
{\rm R^2}$ whose sides are contained in $W^s(\widetilde{Q})$ and $W^u( 
\widetilde{Q}),$ see (\ref{defreten2l}). In order to construct the sets $%
\theta _{\overrightarrow{v}},$ let us first consider the following basic
pieces, denoted $\Gamma _{(1,0)}$ and $\Gamma _{(0,1)}$ (suppose some $
\widetilde{Q}\in p^{-1}(Q)$ is fixed):

\begin{enumerate}
\item  {\bf $\Gamma _{(1,0)}$} is given by the union of $D_{\widetilde{Q}}$
with $D_{\widetilde{Q}}+(1,0)$ and the region bounded by them and two simple
arcs $\eta _1^H,\eta _2^H$ defined as follows: the arc $\eta _1^H$ starts at 
$\widetilde{Q},$ goes through $\widetilde{f}^{N^{\prime \prime }}(\beta _{
\widetilde{Q}})$ until it crosses $\alpha _{\widetilde{Q}}+(1,0)$ and $%
\gamma _{\widetilde{Q}}+(1,0).$ The arc $\eta _2^H$ starts at $\widetilde{f}%
^{N^{\prime \prime }}(\delta _{\widetilde{Q}}\cap \alpha _{\widetilde{Q}}),$
goes through $\widetilde{f}^{N^{\prime \prime }}(\delta _{\widetilde{Q}})$
until it crosses $\alpha _{\widetilde{Q}}+(1,0)$ and $\gamma _{\widetilde{Q}%
}+(1,0).$ We say that the beginning of $\Gamma _{(1,0)}$ is at $D_{
\widetilde{Q}}$ and the end is at $D_{\widetilde{Q}}+(1,0).$

\item  $\Gamma _{(0,1)}$ is given by the union of $D_{\widetilde{Q}}$ with $%
D_{\widetilde{Q}}+(0,1)$ and the region bounded by them and two simple arcs $%
\eta _1^V,\eta _2^V$ analogously defined: the arc $\eta _1^V$ starts at $
\widetilde{Q},$ goes through $\widetilde{f}^{N^{\prime \prime }}(\beta _{
\widetilde{Q}})$ until it crosses $\alpha _{\widetilde{Q}}+(0,1)$ and $%
\gamma _{\widetilde{Q}}+(0,1).$ The arc $\eta _2^V$ starts at $\widetilde{f}%
^{N^{\prime \prime }}(\delta _{\widetilde{Q}}\cap \alpha _{\widetilde{Q}}),$
goes through $\widetilde{f}^{N^{\prime \prime }}(\delta _{\widetilde{Q}})$
until it crosses $\alpha _{\widetilde{Q}}+(0,1)$ and $\gamma _{\widetilde{Q}%
}+(0,1),$ see figure 5. As above, we say that the beginning of $\Gamma
_{(0,1)}$ is at $D_{\widetilde{Q}}$ and the end is at $D_{\widetilde{Q}%
}+(0,1).$
\end{enumerate}

Remember that the definition of $N^{\prime \prime }$ appears in expression (%
\ref{defN2lin}) and below it. Also note that all crosses mentioned above are
topologically transverse intersections in the sense of theorem \ref
{casogeral}.

\begin{lemma}
\label{supline}: Given a vector $\overrightarrow{v}\in {\rm I}%
\negthinspace 
{\rm R^2,}$ we can construct a path connected closed set $\theta _{
\overrightarrow{v}}\subset {\rm I}\negthinspace {\rm R^2}$ such that $\theta
_{\overrightarrow{v}}$ is obtained by the union of integer translates of $%
\Gamma _{(1,0)}$ and $\Gamma _{(0,1)}$ in a way that:\ 
\end{lemma}

\begin{enumerate}
\item  $\theta _{\overrightarrow{v}}$ intersects every straight line
parallel to $\overrightarrow{v}^{\perp },$ a vector orthogonal to $
\overrightarrow{v};$

\item  $\theta _{\overrightarrow{v}}$ is bounded in the direction of $
\overrightarrow{v}^{\perp },$ that is, $\theta _{\overrightarrow{v}}$ is
contained between two straight lines $l_{-}$ and $l_{+},$ both parallel to $
\overrightarrow{v},$ and the distance between these lines is less then $%
3+2.\max \{diameter(\Gamma _{(1,0)}),diameter(\Gamma _{(0,1)})\}.$ So, in
particular $(\theta _{\overrightarrow{v}})^c$ has at least two unbounded
connected components, one containing $l_{-}$ and the other containing $l_{+};
$
\end{enumerate}

{\it Proof:}

To prove this lemma, we fix some $\widetilde{Q}\in {\rm I}\negthinspace {\rm %
R^2}$ as in theorem \ref{casogeral} and consider a straight line $r$ passing
through $\widetilde{Q}$ parallel to $\overrightarrow{v}.$ Without loss of
generality, we can assume that $\widetilde{Q}=(0,0)$ and $\overrightarrow{v}%
=(a,b)$ (so let $\overrightarrow{v}^{\perp }=(-b,a)),$ with $a\geq 0,b\in 
{\rm I}\negthinspace {\rm R}$ and $a^2+b^2=1.$ If $a=0,$ then 
$$
\begin{array}{c}
\theta _{
\overrightarrow{v}}=\stackunder{i\in integers}{\cup }\left( \Gamma
_{(0,1)}+(0,i)\right) \\ \text{and if }b=0, \\ \theta _{\overrightarrow{v}}=%
\stackunder{i\in integers}{\cup }\left( \Gamma _{(1,0)}+(i,0)\right) , 
\end{array}
$$
so first, let us consider the case $a,b>0.$ We denote the Euclidean distance
between two points in the plane by $d_{Euc}(\bullet ,\bullet ).$

We start building the piece of $\theta _{\overrightarrow{v}}$ which follows
the semi-line contained in $r$ given by $\left\{ y=(b/a).x:x\geq 0\right\} .$
Our strategy is the following. We compute the numbers

\begin{equation}
\label{distzero} 
\begin{array}{c}
\left| \left\langle 
\overrightarrow{v}^{\perp },(1,0)\right\rangle \right| \stackrel{def.}{=}%
a_0=d_{Euc}(\widetilde{Q}+(1,0),r)=\left| -b\right| \\ \text{and} \\ \left|
\left\langle \overrightarrow{v}^{\perp },(0,1)\right\rangle \right| 
\stackrel{def.}{=}b_0=d_{Euc}(\widetilde{Q}+(0,1),r)=a. 
\end{array}
\end{equation}
If $a_0\leq b_0,$ then we start with $\Gamma _{(1,0)}.$ In this case $n_0%
\stackrel{def.}{=}(1,0).$ If $a_0>b_0,$ then we start with $\Gamma _{(0,1)}.$
In this case $n_0\stackrel{def.}{=}(0,1).$

So we have our first approximation, which is $\theta _{\overrightarrow{v}%
}^{0+}\stackrel{def.}{=}\Gamma _{n_0},$ where $n_0\in \{(0,1),(1,0)\}$ is
chosen as explained above. Now, in order to decide which of the subsets, $%
\Gamma _{(1,0)}+n_0$ or $\Gamma _{(0,1)}+n_0$ we add, we make the following
computations analogous to the ones in (\ref{distzero}):%
$$
\begin{array}{c}
\left| \left\langle 
\overrightarrow{v}^{\perp },n_0+(1,0)\right\rangle \right| \stackrel{def.}{=}%
a_1=d_{Euc}(\widetilde{Q}+n_0+(1,0),r) \\ \text{and} \\ \left| \left\langle 
\overrightarrow{v}^{\perp },n_0+(0,1)\right\rangle \right| \stackrel{def.}{=}%
b_1=d_{Euc}(\widetilde{Q}+n_0+(0,1),r). 
\end{array}
$$

If $a_1\leq b_1,$ then we add $\Gamma _{(1,0)}+n_0.$ If $a_1>b_1,$ then we
add $\Gamma _{(0,1)}+n_0.$ Now we have $\theta _{\overrightarrow{v}}^{1+}%
\stackrel{def.}{=}\Gamma _{n_0}\cup (\Gamma _{n_1}+n_0),$ where as before $%
n_1\in \{(0,1),(1,0)\}.$ Continuing, in order to decide which of the subsets 
$\Gamma _{(1,0)}+n_0+n_1$ or $\Gamma _{(0,1)}+n_0+n_1$ we add, we compute:%
$$
\begin{array}{c}
\left| \left\langle 
\overrightarrow{v}^{\perp },n_0+n_1+(1,0)\right\rangle \right| \stackrel{def.%
}{=}a_2=d_{Euc}(\widetilde{Q}+n_0+n_1+(1,0),r) \\ \text{and} \\ \left|
\left\langle \overrightarrow{v}^{\perp },n_0+n_1+(0,1)\right\rangle \right| 
\stackrel{def.}{=}b_2=d_{Euc}(\widetilde{Q}+n_0+n_1+(0,1),r). 
\end{array}
$$

If $a_2\leq b_2,$ then we add $\Gamma _{(1,0)}+n_0+n_1.$ If $a_2>b_2,$ then
we add $\Gamma _{(0,1)}+n_0+n_1.$ Now we have $\theta _{\overrightarrow{v}%
}^{2+}\stackrel{def.}{=}\Gamma _{n_0}\cup (\Gamma _{n_1}+n_0)\cup (\Gamma
_{n_2}+n_0+n_1),$ again for some $n_2\in \{(0,1),(1,0)\}.$ After $l$ steps
we arrive at 
$$
\theta _{\overrightarrow{v}}^{l+}\stackrel{def.}{=}\Gamma _{n_0}\cup (\Gamma
_{n_1}+n_0)\cup ...\cup (\Gamma _{n_l}+n_0+n_1+...+n_{l-1}). 
$$

By construction, the points $\widetilde{Q},\widetilde{Q}+n_0,\widetilde{Q}%
+n_0+n_1,...,\widetilde{Q}+n_0+n_1+...+n_l$ all belong to $\theta _{
\overrightarrow{v}}^{l+}.$ Now let us prove that for all integers $l\geq 0,$ 
$d_{Euc}(\widetilde{Q}+n_0+n_1+...+n_l,r)\leq 1.$ Clearly, $d_{Euc}( 
\widetilde{Q},r)=0$ and $d_{Euc}(\widetilde{Q}+n_0,r)\leq \min \{\left|
a\right| ,\left| b\right| \}\leq 1.$ So, suppose by induction that for some
integer $i^{\prime }\geq 0,$ $d_{Euc}(\widetilde{Q}+n_0+n_1+...+n_i,r)\leq
1, $ for all $0\leq i\leq i^{\prime }.$ This means that if we define $\Delta
_{i^{\prime }}\stackrel{def.}{=}\left\langle \overrightarrow{v}^{\perp
},n_0+n_1+...+n_{i^{\prime }}\right\rangle ,$ then $\left| \Delta
_{i^{\prime }}\right| \leq 1.$

If $\Delta _{i^{\prime }}>0,$ then 
$$
-1\leq \left\langle \overrightarrow{v}^{\perp },n_0+n_1+...+n_{i^{\prime
}}+(1,0)\right\rangle =\Delta _{i^{\prime }}-b<\Delta _{i^{\prime }}\leq 1. 
$$

If $\Delta _{i^{\prime }}<0,$ then 
$$
-1\leq \Delta _{i^{\prime }}<\left\langle \overrightarrow{v}^{\perp
},n_0+n_1+...+n_{i^{\prime }}+(0,1)\right\rangle =\Delta _{i^{\prime
}}+a\leq 1. 
$$
These estimates clearly imply that $d_{Euc}(\widetilde{Q}+n_0+n_1+...+n_{i^{%
\prime }+1},r)\leq 1$ because $n_{i^{\prime }+1}\in \{(0,1),(1,0)\}$ is
chosen in a way to minimize the distance. So, our claim is proved.

If $\Delta _{i^{\prime }}=0,$ this means that $\widetilde{Q}%
+n_0+n_1+...+n_{i^{\prime }}$ belongs to $r,$ which means that $
\overrightarrow{v}$ is a rational direction and so%
$$
\theta _{\overrightarrow{v}}=\stackunder{i\in integers}{\cup }\left( 
\begin{array}{c}
\Gamma _{n_0}\cup (\Gamma _{n_1}+n_0)\cup ...\cup (\Gamma _{n_{i^{\prime
}}}+n_0+n_1+...+n_{i^{\prime }-1})+ \\ 
+i.(n_0+n_1+...+n_{i^{\prime }-1}+n_{i^{\prime }}) 
\end{array}
\right) 
$$

In case $\Delta _i\neq 0$ for all integers $i>0,$ we define $\theta _{
\overrightarrow{v}}^{+}\stackrel{def.}{=}\stackunder{i\geq 0}{\cup }\theta _{
\overrightarrow{v}}^{i+}.$ In order to get the whole $\theta _{
\overrightarrow{v}},$ we have to construct the other side of it. For this,
let

$$
\theta _{\overrightarrow{v}}^{i-}\stackrel{def.}{=}(\Gamma _{n_0}-n_0)\cup
(\Gamma _{n_1}-n_1-n_0)\cup ...\cup (\Gamma _{n_i}-n_i-n_{i-1}...-n_1-n_0) 
$$
and analogously $\theta _{\overrightarrow{v}}^{-}\stackrel{def.}{=}%
\stackunder{i\geq 0}{\cup }\theta _{\overrightarrow{v}}^{i-}.$ As we did
above, for any integer $i\geq 0$ points of the form $\widetilde{Q}%
-n_0-n_1...-n_i$ all belong to $\theta _{\overrightarrow{v}}^{-}$ and 
$$
d_{Euc}(\widetilde{Q}-n_0-n_1...-n_i,r)=\left| \left\langle \overrightarrow{v%
}^{\perp },-n_0-n_1...-n_i\right\rangle \right| = 
$$

$$
=\left| \left\langle \overrightarrow{v}^{\perp
},n_0+n_1+...+n_i\right\rangle \right| =\left| \Delta _i\right| \leq 1. 
$$

So, finally we make $\theta _{\overrightarrow{v}}\stackrel{def.}{=}\theta _{
\overrightarrow{v}}^{-}\cup \theta _{\overrightarrow{v}}^{+}.$ It is a
closed, connected subset of the plane and from the properties obtained
above, the projection of $\theta _{\overrightarrow{v}}$ in the direction of $
\overrightarrow{v}^{\perp }$ has diameter smaller than $3+2.\max
\{diameter(\Gamma _{(0,1)}),diameter(\Gamma _{(1,0)})\},$ so it is contained
between two straight lines parallel to $\overrightarrow{v},$ whose distance
is less than $3+2.\max \{diameter(\Gamma _{(0,1)}),diameter(\Gamma
_{(1,0)})\}.$ The fact that $\theta _{\overrightarrow{v}}$ intersects every
straight line parallel to $\overrightarrow{v}^{\perp }$ is easy. If $a>0$
and $b<0,$ the proof is analogous. $\Box $

\vskip 0.2truecm

The next lemma uses theorem \ref{main1} and easily implies theorem \ref
{main2}:

\begin{lemma}
\label{fabio}: Suppose $f\in Diff_0^{1+\epsilon }({\rm T^2})$ has a rotation
set $\rho (\widetilde{f})$ with interior. Let $\mu \in M_{inv}(f)$ be such
that the rotation vector of $\mu ,$ $\rho (\mu )\in \partial \rho (
\widetilde{f}).$ Let $r$ be a supporting line at $\rho (\mu )$ and $
\overrightarrow{v}^{\perp }$ be the unitary vector orthogonal to $r,$
pointing towards the connected component of $r^c$ which does not intersect $%
\rho (\widetilde{f}).$ Then, if $x^{\prime }\in supp(\mu ),$ for any $
\widetilde{x}^{\prime }\in p^{-1}(x^{\prime })$ and any integer $n>0,$

\begin{equation}
\label{fff}\left| \left\langle \widetilde{f}^n(\widetilde{x}^{\prime })-
\widetilde{x}^{\prime }-n.\rho (\mu ),\overrightarrow{v}^{\perp
}\right\rangle \right| \leq 2+M_f,
\end{equation}
where $M_f$ comes from theorem \ref{main1}.
\end{lemma}

{\it Proof:}

Let us denote $r^c=\Omega _1\cup \Omega _2,$ in a way that $\rho (\widetilde{%
f})\subset r\cup \Omega _1.$

\begin{description}
\item[Fact]  \label{ergdecomp} {\bf \ref{ergdecomp}}: Every ergodic measure $%
\xi $ that appears in the ergodic decomposition of $\mu $ has rotation
vector contained in $r.$
\end{description}

{\it Proof:}

This follows from $\rho (\widetilde{f})\cap \Omega _2=\emptyset $ and $\rho
(\mu )\in r.$ By contradiction, assume that for some $\xi $ in the ergodic
decomposition of $\mu ,$ $\rho (\xi )$ does not belong to $r.$ Then $\rho
(\xi )\in \Omega _1.$ Here we are using the non-obvious fact that%
$$
\rho (\widetilde{f})=\{\omega \in {\rm I\negthinspace R^2:}\exists \eta \in
M_{inv}(f)\text{ such that }\rho (\eta )=\int_{{\rm T^2}}\phi (x)d\eta
=\omega \}, 
$$
see \cite{misiu}. Therefore, as $\xi $ is in the ergodic decomposition of $%
\mu ,$ the fact that $\rho (\xi )\in \Omega _1$ would imply the existence of
another ergodic measure $\xi ^{\prime }$ also in the ergodic decomposition
of $\mu $ such that $\rho (\xi ^{\prime })\in \Omega _2$ (because $\rho (\mu
)\in r).$ This contradiction proves the fact. $\Box $

\vskip 0.2truecm

To prove lemma \ref{fabio}, we again argue by contradiction. So let us
suppose that there exists $x^{\prime }\in supp(\mu )$ and some integer $%
n_0>0,$ such that for any $\widetilde{x}^{\prime }\in p^{-1}(x^{\prime }),$

\begin{equation}
\label{contrafff}\left\langle \widetilde{f}^{n_0}(\widetilde{x}^{\prime })- 
\widetilde{x}^{\prime }-n_0.\rho (\mu ),\overrightarrow{v}^{\perp
}\right\rangle <-2-M_f. 
\end{equation}

Theorem \ref{main1} implies that if the present lemma does not hold, then
the above is the only possibility.

Expression (\ref{contrafff}) and a simple continuity argument clearly imply
that there exists $\epsilon ^{\prime }>0$ such that for all $x\in
B_{\epsilon ^{\prime }}(x^{\prime })$ (the ball of radius $\epsilon ^{\prime
}$ centered at $x^{\prime })$ and any $\widetilde{x}\in p^{-1}(x),$%
\begin{equation}
\label{contrafffbola}\left\langle \widetilde{f}^{n_0}(\widetilde{x})- 
\widetilde{x}-n_0.\rho (\mu ),\overrightarrow{v}^{\perp }\right\rangle
<-2-M_f. 
\end{equation}

Now let $\nu \in M_{inv}(f)$ be an ergodic measure in the ergodic
decomposition of $\mu $ such that $x^{\prime }\in supp(\nu ).$ As $\rho (\nu
)\in r$ (see fact \ref{ergdecomp}), $\rho (\nu )=\rho (\mu )+\lambda . 
\overrightarrow{v},$ where $\overrightarrow{v}$ is parallel to $r$ and $%
\lambda $ is some adequate real number. So, $\left\langle \rho (\nu ), 
\overrightarrow{v}^{\perp }\right\rangle =\left\langle \rho (\mu ), 
\overrightarrow{v}^{\perp }\right\rangle .$

We also define the relative to $\mu $ displacement function in the direction
of $\overrightarrow{v}^{\perp }$ as $\phi _{\mu ,\overrightarrow{v}^{\perp
}}:{\rm T^2\rightarrow I\negthinspace R}$ given by $\phi _{\mu , 
\overrightarrow{v}^{\perp }}(x)=\left\langle \widetilde{f}(\widetilde{x})- 
\widetilde{x}-\rho (\mu ),\overrightarrow{v}^{\perp }\right\rangle ,$ for
any $\widetilde{x}\in p^{-1}(x).$ Then the following consequences hold:

\begin{enumerate}
\item  $\int_{{\rm T^2}}\phi _{\mu ,\overrightarrow{v}^{\perp }}(x)d\nu =0;$

\item  for any $\widetilde{x}\in {\rm I\negthinspace R^2}$ and any integer $%
n>0,$ if $x=p(\widetilde{x}),$ then

$$
\left\langle \widetilde{f}^n(\widetilde{x})-\widetilde{x}-n.\rho (\mu ),
\overrightarrow{v}^{\perp }\right\rangle =\stackrel{n-1}{\stackunder{i=0}{%
\sum }}\phi _{\mu ,\overrightarrow{v}^{\perp }}(f^i(x)); 
$$
\end{enumerate}

So from Atkinson's lemma (see \cite{atkinson}) we get that for every $%
0<\epsilon <\epsilon ^{\prime },$ there exists $x^{*}\in B_\epsilon
(x^{\prime }),$ such that for some integer $n_1>n_0$ and any $\widetilde{x}%
^{*}\in p^{-1}(x^{*}),$%
$$
\left| \left\langle \widetilde{f}^{n_1}(\widetilde{x}^{*})-\widetilde{x}%
^{*}-n_1.\rho (\mu ),\overrightarrow{v}^{\perp }\right\rangle \right| <1. 
$$
Thus, from expressions (\ref{contrafffbola}) and the above one, we finally
obtain that 
$$
\left\langle \widetilde{f}^{n_1-n_0}(\widetilde{f}^{n_0}(\widetilde{x}%
^{*}))- \widetilde{f}^{n_0}(\widetilde{x}^{*})-(n_1-n_0).\rho (\mu ), 
\overrightarrow{v}^{\perp }\right\rangle >1+M_f, 
$$
a contradiction with theorem \ref{main1}. So expression (\ref{contrafff})
does not hold and the lemma is proved. $\Box $

\vskip 0.2truecm

\subsection{Proof of theorem 1}

First, let us consider a map $\widetilde{g}(\bullet )\stackrel{def.}{=} 
\widetilde{f}^q(\bullet )-(p,s)$ for some rational vector $\left( \frac
pq,\frac sq\right) \in int(\rho (\widetilde{f})),$ not necessarily in
irreducible form, in a way that $g$ has a FIXED hyperbolic saddle point $%
Q\in {\rm T^2}$ with positive eigenvalues, as in theorem \ref{casogeral}.
For example, $\left( \frac pq,\frac sq\right) $ could be equal to $\left(
\frac 13,\frac 23\right) ,$ but $q=30,p=10$ and $s=20.$ It is easy to see
that $\rho (\widetilde{g})=q.\rho (\widetilde{f})-(p,s).$ So if we fix some $%
\omega \in \partial \rho (\widetilde{f})$ and a supporting line $r$ at $%
\omega ,$ parallel to some unitary vector $\overrightarrow{v},$ the
corresponding rotation vector and supporting line for $\widetilde{g}$ are: $%
q.\omega -(p,s)\in \partial \rho (\widetilde{g})$ and a straight line $%
r^{\prime }$ passing through $q.\omega -(p,s),$ also parallel to $
\overrightarrow{v}.$

Let us show that if the theorem holds for $g,$ then it also holds for $f.$
For this, assume there exists a number $M_g>0$ such that for any $\tau \in
\partial \rho (\widetilde{g})$ and any supporting line $r$ at $\tau $$,$ if $
\overrightarrow{v}^{\perp }$ is the unitary vector orthogonal to $r,$
pointing towards the connected component of $r^c$ which does not intersect $%
\rho (\widetilde{g}),$ then

\begin{equation}
\label{teo1g}\left\langle \widetilde{g}^n(\widetilde{x})-\widetilde{x}%
-n.\tau ,\overrightarrow{v}^{\perp }\right\rangle \leq M_g,\text{ for all } 
\widetilde{x}\in {\rm I\negthinspace R^2}\text{ and any integer }n>0. 
\end{equation}

From the relation between $\rho (\widetilde{g})$ and $\rho (\widetilde{f}),$ 
$$
\omega \in \partial \rho (\widetilde{f})\text{ }\Leftrightarrow \text{ }%
q.\omega -(p,s)\in \partial \rho (\widetilde{g}). 
$$

Expression (\ref{teo1g}) implies that 
$$
\left\langle \widetilde{f}^{n.q}(\widetilde{x})-\widetilde{x}-n.q.\frac{\tau
+(p,s)}q,\overrightarrow{v}^{\perp }\right\rangle \leq M_g,\text{ for all } 
\widetilde{x}\in {\rm I\negthinspace R^2}\text{ and any integer }n>0. 
$$

Which gives,

$$
\begin{array}{c}
\left\langle 
\widetilde{f}^n(\widetilde{x})-\widetilde{x}-n.\frac{\tau +(p,s)}q, 
\overrightarrow{v}^{\perp }\right\rangle \leq M_g+q.\left( \stackunder{ 
\widetilde{z}\in {\rm I\negthinspace R^2}}{\sup }\left\| \widetilde{f}( 
\widetilde{z})-\widetilde{z}\right\| +\stackunder{\iota \in \rho (\widetilde{%
f})}{\sup }\left\| \iota \right\| \right) , \\  \\ 
\text{ for all }\widetilde{x}\in {\rm I\negthinspace R^2}\text{ and any
integer }n>0. 
\end{array}
$$

As $\stackunder{\iota \in \rho (\widetilde{f})}{\sup }\left\| \iota \right\|
\leq \stackunder{\widetilde{z}\in {\rm I\negthinspace R^2}}{\sup }\left\| 
\widetilde{f}(\widetilde{z})-\widetilde{z}\right\| $ and the map $\tau $$%
\rightarrow \frac{\tau +(p,s)}q$ is a bijection from $\partial \rho (
\widetilde{g})$ to $\partial $$\rho (\widetilde{f}),$ if we choose $%
M_f=M_g+2q.\left( \stackunder{\widetilde{z}\in {\rm I\negthinspace R^2}}{%
\sup }\left\| \widetilde{f}(\widetilde{z})-\widetilde{z}\right\| \right) ,$
then we are done.

So it remains to show that the present theorem holds for $g.$ Let us fix
some $\tau \in \partial \rho (\widetilde{g})$ and any supporting line $r$ at 
$\tau $$,$ parallel to some vector $\overrightarrow{v}.$ Also, let $
\overrightarrow{v}^{\perp }$ be the unitary vector orthogonal to $r,$
pointing towards the connected component of $r^c$ which does not intersect $%
\rho (\widetilde{g}).$ From lemma \ref{supline}, fixed some $\widetilde{Q}%
\in p^{-1}(Q),$ there exists a subset $\theta _{\overrightarrow{v}}\subset 
{\rm I}\negthinspace {\rm R^2}$ as in the statement of that lemma,
containing $\widetilde{Q}.$ This means that there are straight lines, $l_{-}$
and $l_{+},$ both parallel to $\overrightarrow{v},$ at a distance less then $%
d_\theta \stackrel{def.}{=}3+2.\max \{diameter(\Gamma
_{(0,1)}),diameter(\Gamma _{(1,0)})\}$ such that $\theta _{\overrightarrow{v}%
}$ is contained between $l_{-}$ and $l_{+}$ and $\overrightarrow{v}^{\perp }$
points from $l_{-}$ to $l_{+}.$ Let $U_{\overrightarrow{v}}^{-}$ and $U_{
\overrightarrow{v}}^{+}$ be the unbounded connected components of $(\theta _{
\overrightarrow{v}})^c,$ such that $\overrightarrow{v}^{\perp }$ points
towards $U_{\overrightarrow{v}}^{+},$ or equivalently $l_{-}\subset U_{
\overrightarrow{v}}^{-}$ and $l_{+}\subset U_{\overrightarrow{v}}^{+}.$ Note
that $\theta _{\overrightarrow{v}}$ also intersects all straight lines
parallel to $\overrightarrow{v}^{\perp }$.

If $(a,b)$ is an integer vector such that $\left| \left\langle (a,b), 
\overrightarrow{v}^{\perp }\right\rangle \right| >d_\theta ,$ then $\theta _{
\overrightarrow{v}}\cap (\theta _{\overrightarrow{v}}+(a,b))=\emptyset .$
More precisely, if $\left\langle (a,b),\overrightarrow{v}^{\perp
}\right\rangle >d_\theta ,$ then $(\theta _{\overrightarrow{v}%
}+(a,b))\subset U_{\overrightarrow{v}}^{+}$ and if $\left\langle (a,b), 
\overrightarrow{v}^{\perp }\right\rangle <-d_\theta ,$ then $(\theta _{
\overrightarrow{v}}+(a,b))\subset U_{\overrightarrow{v}}^{-}.$

Now let us suppose by contradiction that there exists $\widetilde{x}^{*}\in 
{\rm I\negthinspace R^2}$ and an integer $n^{*}>N^{\prime \prime }>0$ such
that 
\begin{equation}
\label{expprin}\left\langle \widetilde{g}^{n^{*}}(\widetilde{x}^{*})- 
\widetilde{x}^{*}-n^{*}.\tau ,\overrightarrow{v}^{\perp }\right\rangle
>100+20.d_\theta +N^{\prime \prime }.\left\langle \tau ,\overrightarrow{v}%
^{\perp }\right\rangle . 
\end{equation}

Remember that $N^{\prime \prime }$ was defined in expression (\ref{defN2lin}%
) and $\left\langle \tau ,\overrightarrow{v}^{\perp }\right\rangle >0$
because $(0,0)\in int(\rho (\widetilde{g}))$. Without loss of generality, we
can assume that $\widetilde{x}^{*}$ belongs to the connected component of $%
(l_{-})^c$ contained in $U_{\overrightarrow{v}}^{-}$ and $d_{Euc}( 
\widetilde{x}^{*},l_{-})<1.$

Now let us choose some integer vector $(a_{+},b_{+})$ such that 
$$
50+3.d_\theta +(n^{*}+N^{\prime \prime }).\left\langle \tau ,\overrightarrow{%
v}^{\perp }\right\rangle <\left\langle (a_{+},b_{+}),\overrightarrow{v}%
^{\perp }\right\rangle <70+10.d_\theta +(n^{*}+N^{\prime \prime
}).\left\langle \tau ,\overrightarrow{v}^{\perp }\right\rangle . 
$$
From what we explained above, $(\theta _{\overrightarrow{v}%
}+(a_{+},b_{+}))\subset U_{\overrightarrow{v}}^{+}$ and $\widetilde{g}%
^{n^{*}}(\widetilde{x}^{*})$ belongs to the connected component of $%
(l_{+}+(a_{+},b_{+}))^c $ contained in $U_{\overrightarrow{v}%
}^{+}+(a_{+},b_{+}).$ So, $\widetilde{g}^{n^{*}}(\theta _{\overrightarrow{v}%
})$ intersects $\theta _{\overrightarrow{v}}+(a_{+},b_{+}).$ More precisely,
there exist integer vectors $(a_i,b_i),(a_f,b_f)$ such that $\widetilde{Q}%
+(a_i,b_i)\in \theta _{\overrightarrow{v}},$ $\widetilde{Q}+(a_f,b_f)\in
\theta _{\overrightarrow{v}}+(a_{+},b_{+})$ and at least one of the
following possibilities hold:

\begin{itemize}
\item  $\widetilde{g}^{n^{*}}(\eta _1^H+(a_i,b_i))$ and $\widetilde{g}%
^{n^{*}}(\eta _2^H+(a_i,b_i))$ intersect both $\alpha _{\widetilde{Q}%
}+(a_f,b_f)$ and $\gamma _{\widetilde{Q}}+(a_f,b_f);$

\item  $\widetilde{g}^{n^{*}}(\eta _1^V+(a_i,b_i))$ and $\widetilde{g}%
^{n^{*}}(\eta _2^V+(a_i,b_i))$ intersect both $\alpha _{\widetilde{Q}%
}+(a_f,b_f)$ and $\gamma _{\widetilde{Q}}+(a_f,b_f);$
\end{itemize}

In any of the above cases, from the definition of $N^{\prime \prime }$ (see
expression (\ref{defN2lin})), we get that $\widetilde{g}^{n^{*}+N^{\prime
\prime }}(D_{\widetilde{Q}}+(a_i,b_i))\cap (D_{\widetilde{Q}}+(a_f,b_f))$
contains a topological rectangle $\widetilde{R}_{fast}$ with one side
contained in $\alpha _{\widetilde{Q}}+(a_f,b_f),$ another one contained in $%
\gamma _{\widetilde{Q}}+(a_f,b_f)$ and the two other sides contained in the
interior of $D_{\widetilde{Q}}+(a_f,b_f),$ see figure 6. This implies that
there is a compact $g$-invariant set contained in $D_Q=p(D_{\widetilde{Q}%
})\subset {\rm T^2}$ whose dynamics is semi-conjugate to that of a
horseshoe. In particular, there exists a fixed point for $g^{n^{*}+N^{\prime
\prime }}$ in $p(\widetilde{R}_{fast})$ (such a point must exist, but it may
not be unique) and its rotation vector with respect to $\widetilde{g}$ is 
$$
\rho _{fast}=\left( \frac{a_f-a_i,b_f-b_i}{n^{*}+N^{\prime \prime }}\right)
. 
$$

As $\widetilde{Q},\widetilde{Q}+(a_i,b_i)\in \theta _{\overrightarrow{v}}$
and $\widetilde{Q}+(a_f,b_f)\in \theta _{\overrightarrow{v}}+(a_{+},b_{+})$
we get that 
$$
\left\langle \left( \frac{a_f-a_i,b_f-b_i}{n^{*}+N^{\prime \prime }}\right)
-\tau ,\overrightarrow{v}^{\perp }\right\rangle =\left\langle \left( \frac{ 
\widetilde{Q}+(a_f,b_f)-(\widetilde{Q}+(a_i,b_i))}{n^{*}+N^{\prime \prime }}%
\right) -\tau ,\overrightarrow{v}^{\perp }\right\rangle > 
$$

$$
>\frac{50+3.d_\theta +(n^{*}+N^{\prime \prime }).\left\langle \tau ,
\overrightarrow{v}^{\perp }\right\rangle -2.d_\theta }{n^{*}+N^{\prime
\prime }}-\left\langle \tau ,\overrightarrow{v}^{\perp }\right\rangle =\frac{%
50+d_\theta }{n^{*}+N^{\prime \prime }}>0. 
$$
And this is a contradiction, because from $\left\langle \rho _{fast}-\tau ,
\overrightarrow{v}^{\perp }\right\rangle >0$ we get that $\rho _{fast}\notin
\rho (\widetilde{g}).$ So expression (\ref{expprin}) does not hold and thus,
for all $\widetilde{x}\in {\rm I\negthinspace R^2}$ and any integer $n>0,$ 
$$
\left\langle \widetilde{g}^n(\widetilde{x})-\widetilde{x}-n.\tau ,
\overrightarrow{v}^{\perp }\right\rangle \leq \max \{\left( 100+20.d_\theta
+N^{\prime \prime }.\left\langle \tau ,\overrightarrow{v}^{\perp
}\right\rangle \right) ,\text{ }M_{g,N^{\prime \prime }}\}, 
$$
where the number 
$$
M_{g,N^{\prime \prime }}\stackrel{def.}{=}2N^{\prime \prime }.\left( 
\stackunder{\widetilde{x}\in {\rm I\negthinspace R^2}}{\sup }\left\| 
\widetilde{g}(\widetilde{x})-\widetilde{x}\right\| \right)  
$$
appears because for any $\widetilde{x}\in {\rm I\negthinspace R^2\ }$and any$%
{\rm \ }0\leq n\leq N^{\prime \prime },$ 
$$
\begin{array}{c}
\left\langle 
\widetilde{g}^n(\widetilde{x})-\widetilde{x}-n.\tau ,\overrightarrow{v}%
^{\perp }\right\rangle \leq \left\| \widetilde{g}^n(\widetilde{x})-
\widetilde{x}\right\| +n.\left\| \tau \right\| \leq  \\  \\ 
\leq n.\left( \stackunder{\widetilde{x}\in {\rm I\negthinspace R^2}}{\sup }%
\left\| \widetilde{g}(\widetilde{x})-\widetilde{x}\right\| \right) +n.%
\stackunder{\iota \in \rho (\widetilde{g})}{\sup }\left\| \iota \right\|
\leq 2n.\left( \stackunder{\widetilde{x}\in {\rm I\negthinspace R^2}}{\sup }%
\left\| \widetilde{g}(\widetilde{x})-\widetilde{x}\right\| \right) \leq
M_{g,N^{\prime \prime }}.
\end{array}
$$

So we can take 
$$
\begin{array}{c}
M_g 
\stackrel{def.}{=}\max \{\left( 100+20.d_\theta +N^{\prime \prime }.%
\stackunder{\widetilde{x}\in {\rm I\negthinspace R^2}}{\sup }\left\| 
\widetilde{g}(\widetilde{x})-\widetilde{x}\right\| \right) ,M_{g,N^{\prime
\prime }}\}, \\  
\end{array}
$$
because $\stackunder{\widetilde{x}\in {\rm I\negthinspace R^2}}{\sup }%
\left\| \widetilde{g}(\widetilde{x})-\widetilde{x}\right\| \geq \stackunder{%
\iota \in \rho (\widetilde{g})}{\sup }\left\| \iota \right\| \geq
\left\langle \tau ,\overrightarrow{v}^{\perp }\right\rangle $ for every $%
\tau \in \partial \rho (\widetilde{g})$ and $\overrightarrow{v}^{\perp }$ an
unitary vector orthogonal to the supporting line at $\tau $ oriented in an
adequate way. $\Box $

\vskip 0.2truecm

\subsection{Proof of theorem 2}

By contradiction, suppose that the rotation vector of the Lebesgue measure,
denoted $\omega ,$ belongs to $\partial $$\rho (\widetilde{f}).$ Let $r$ be
a supporting line at $\omega $ and let $\overrightarrow{v}^{\perp }$ be a
unitary vector orthogonal to $r$, pointing towards the connected component
of $r^c$ that does not intersect $\rho (\widetilde{f}).$ From lemma \ref
{fabio} we get that for all $\widetilde{x}\in {\rm I\negthinspace R^2}$ and
any integer $n>0$ (remember that $supp(Lebesgue)={\rm T^2}),$

\begin{equation}
\label{denovo}\left| \left\langle \widetilde{f}^n(\widetilde{x})-\widetilde{x%
}-n.\omega ,\overrightarrow{v}^{\perp }\right\rangle \right| \leq 2+M_f, 
\text{ where }M_f\text{ comes from theorem \ref{main1}.} 
\end{equation}

Now pick some point $z\in {\rm T^2}$ which is periodic and has a rotation
vector $\nu \in int(\rho (\widetilde{f})).$ As $\nu \notin r,$ $\left\langle
\nu ,\overrightarrow{v}^{\perp }\right\rangle \neq \left\langle \omega , 
\overrightarrow{v}^{\perp }\right\rangle .$ So, for any $\widetilde{z}\in
p^{-1}(z),$%
$$
\frac{\left\langle \widetilde{f}^n(\widetilde{z})-\widetilde{z}-n.\omega , 
\overrightarrow{v}^{\perp }\right\rangle }n\stackrel{n\rightarrow \infty }{%
\rightarrow }\left\langle \nu -\omega ,\overrightarrow{v}^{\perp
}\right\rangle \neq 0, 
$$
a contradiction with expression (\ref{denovo}). This proves the theorem. $%
\Box $

\vskip 0.2truecm

\subsection{Proof of corollary 3}

First, assume $\rho (\mu )\in \partial \rho (\widetilde{f})$ is a vertex.
Then there are 2 different supporting lines at $\rho (\mu )$ (in fact, there
are infinitely many), denoted $r_1$ and $r_2,$ and $\overrightarrow{v}%
_1^{\perp },\overrightarrow{v}_2^{\perp }$ are unitary vectors orthogonal,
respectively to $r_1$ and $r_2,$ such that $-\overrightarrow{v}_1^{\perp }$
and $-\overrightarrow{v}_2^{\perp }$ point towards $\rho (\widetilde{f}).$
From lemma \ref{fabio}, for every $x\in supp(\mu ),$ for any $\widetilde{x}%
\in p^{-1}(x)$ and any integer $n>0,$

$$
\left| \left\langle \widetilde{f}^n(\widetilde{x})-\widetilde{x}-n.\rho (\mu
),\overrightarrow{v}_i^{\perp }\right\rangle \right| \leq 2+M_f,\text{ for }%
i=1,2. 
$$

As $\overrightarrow{v}_1^{\perp }$ and $\overrightarrow{v}_2^{\perp }$ are
not parallel, the above expression implies that there exists $M_\mu $ which
depends only on $M_f$ and $\overrightarrow{v}_1^{\perp },\overrightarrow{v}%
_2^{\perp }$ such that for every $x\in supp(\mu ),$ for any $\widetilde{x}%
\in p^{-1}(x)$ and any integer $n>0,$ 
$$
\left\| \widetilde{f}^n(\widetilde{x})-\widetilde{x}-n.\rho (\mu )\right\|
<M_\mu . 
$$
This proves the first part of the corollary. Now suppose $\rho (\mu )$ is an
extremal point of $\rho (\widetilde{f})$ and the intersection of the
supporting line $r$ at $\rho (\mu )$ with $\rho (\widetilde{f})$ is just $%
\rho (\mu ).$ From lemma \ref{fabio}, we know that for every $x\in supp(\mu
) $ and for any $\widetilde{x}\in p^{-1}(x),$ 
\begin{equation}
\label{umamais}\stackunder{n\rightarrow \infty }{\lim }\left\langle \frac{ 
\widetilde{f}^n(\widetilde{x})-\widetilde{x}}n-\rho (\mu ),\overrightarrow{v}%
^{\perp }\right\rangle =0, 
\end{equation}
where $\overrightarrow{v}^{\perp }$ is the unitary vector orthogonal to $r$
oriented in a way that $-\overrightarrow{v}^{\perp }$ points towards $\rho ( 
\widetilde{f}).$ As the accumulation points of the sequence 
$$
\frac{\widetilde{f}^n(\widetilde{x})-\widetilde{x}}n 
$$
belong both to $r$ (this follows from expression (\ref{umamais})) and to $%
\rho (\widetilde{f}),$ there is just one accumulation point and it is $\rho
(\mu ).$ So $\rho (x)$ exists and it is equal to $\rho (\mu ).$ As $x$ is
any point in $supp(\mu ),$ this proves the second part of the corollary. $%
\Box $

\vskip 0.2truecm

{\it Acknowledgements: }I would like to thank F\'abio Tal for conversations
concerning lemma \ref{fabio}.

\centerline{\bf Figure captions.}

\begin{itemize}
\item[Figure 1. ]  Diagram showing a topologically transverse intersection
between $W^u(\widetilde{Q})$ and $W^s(\widetilde{Q}+(a,b)).$

\item[Figure 2.]  Diagram showing the topological rectangle $D_{\widetilde{Q}%
}.$

\item[Figure 3.]  Diagram showing the set $\theta .$

\item[Figure 4.]  Diagram showing the topological rectangle $\widetilde{R}%
^{*}$ contained in $\widetilde{f}^{n_M+N^{\prime \prime }}(D_{\widetilde{Q}%
})\cap \left( D_{\widetilde{Q}}+(a,n_M+\left\lfloor
M-5-4.d_{(1,0)}\right\rfloor )\right) .$

\item[Figure 5.]  Diagram showing the sets: (a) $\Gamma _{(0,1)}$ and (b) $%
\Gamma _{(1,0)}.$

\item[Figure 6.]  Diagram showing the the topological rectangle $\widetilde{R%
}_{fast}.$
\end{itemize}

\begin{center}
\mbox{\includegraphics[width=13cm]{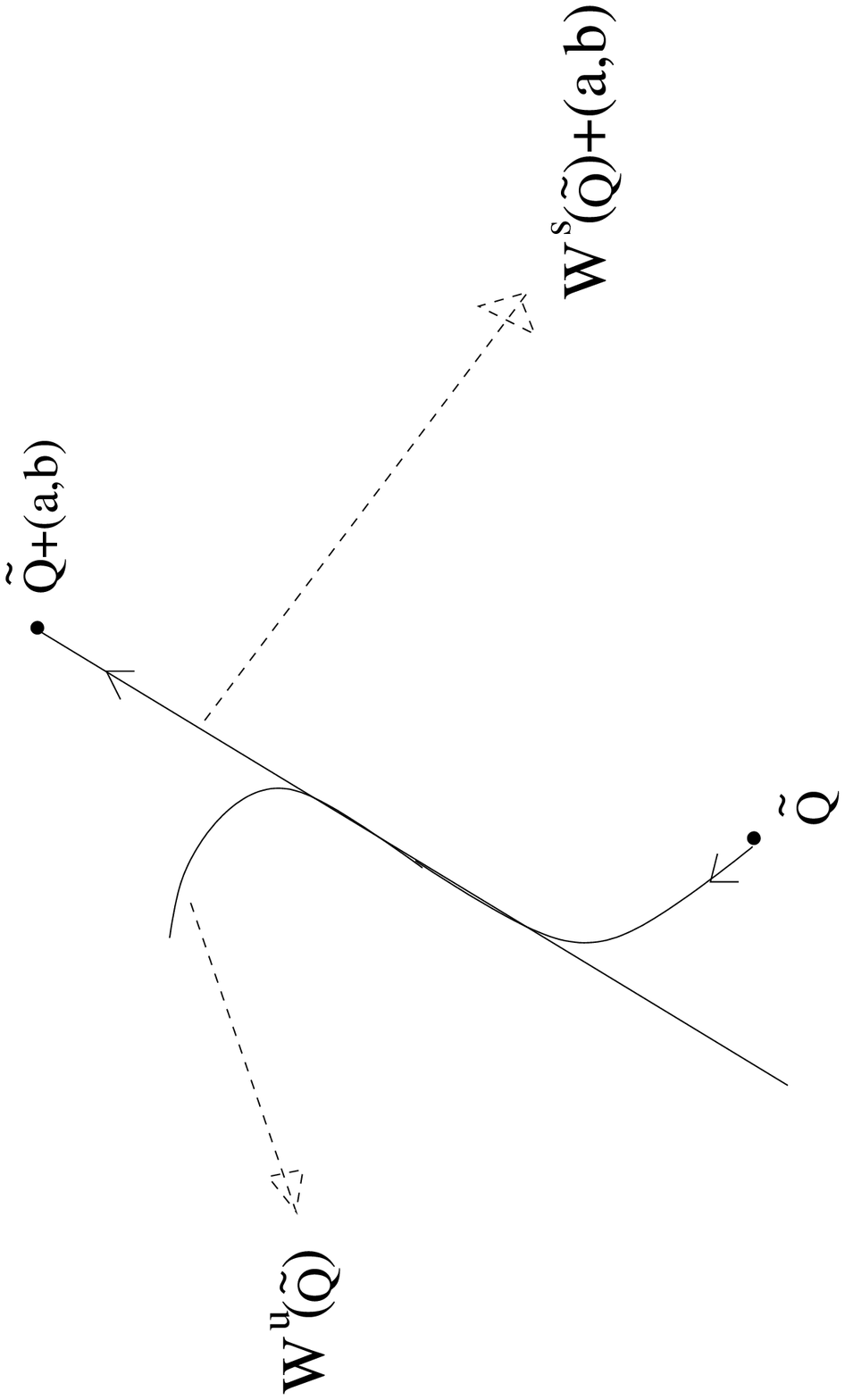}} 
\mbox{\includegraphics[width=13cm]{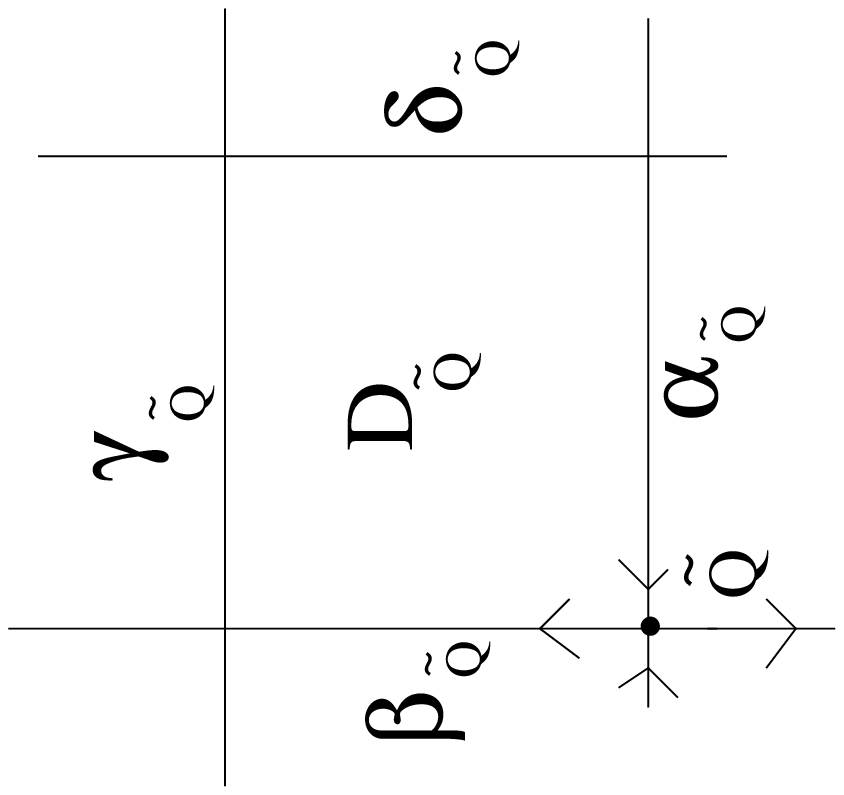}} 
\mbox{\includegraphics[width=13cm]{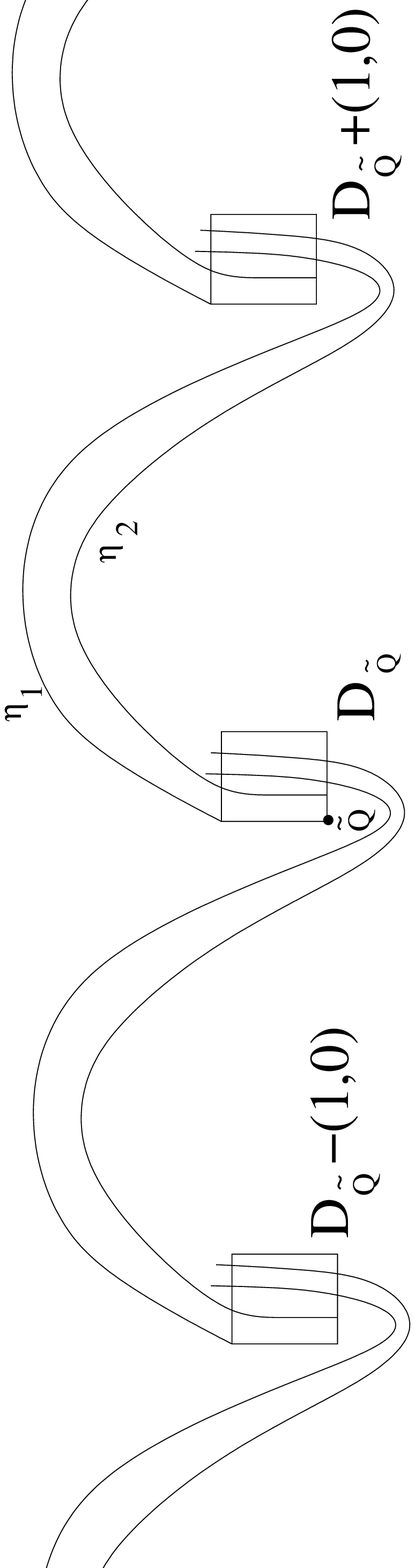}} 
\mbox{\includegraphics[width=13cm]{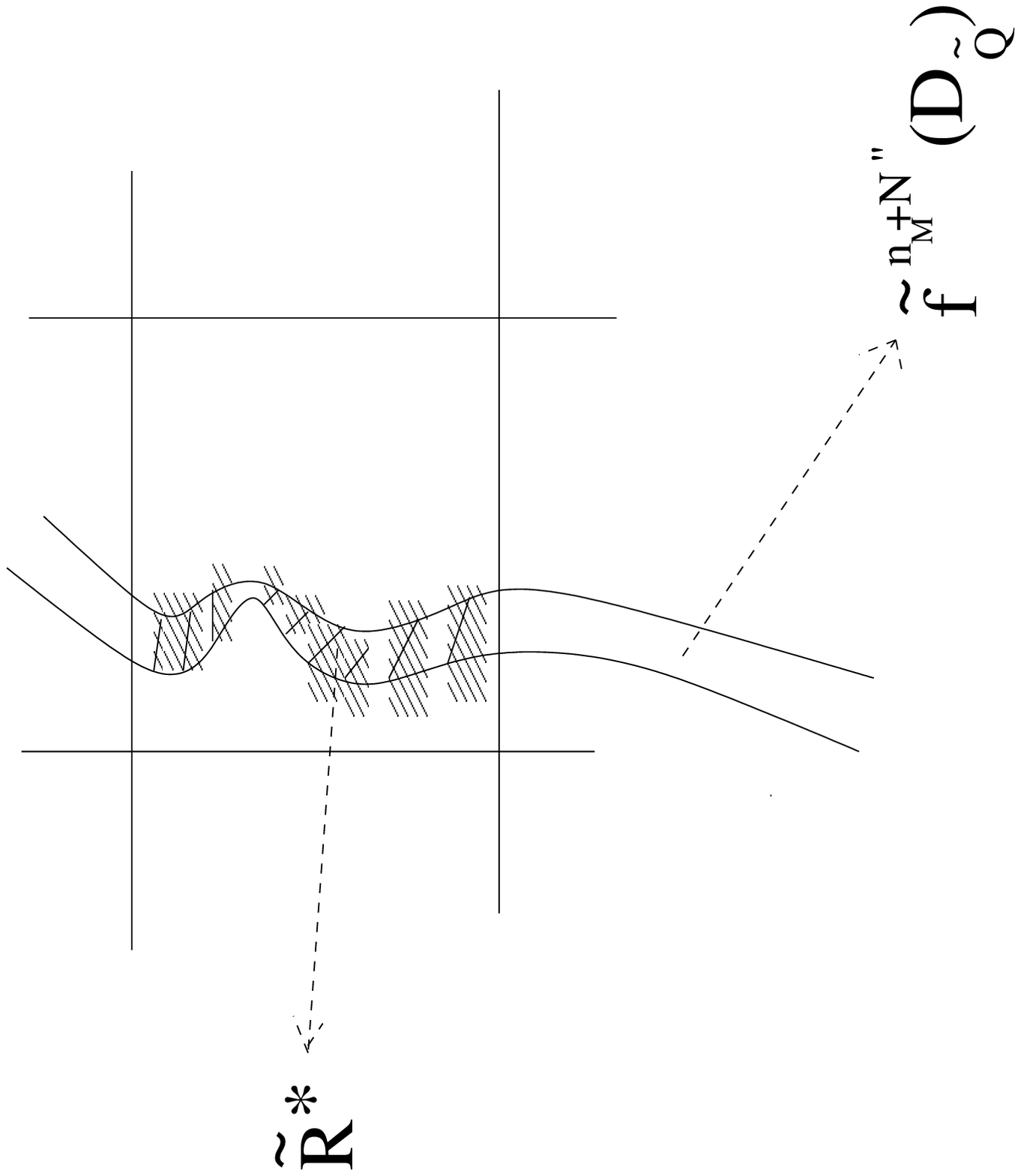}} 
\mbox{\includegraphics[width=13cm]{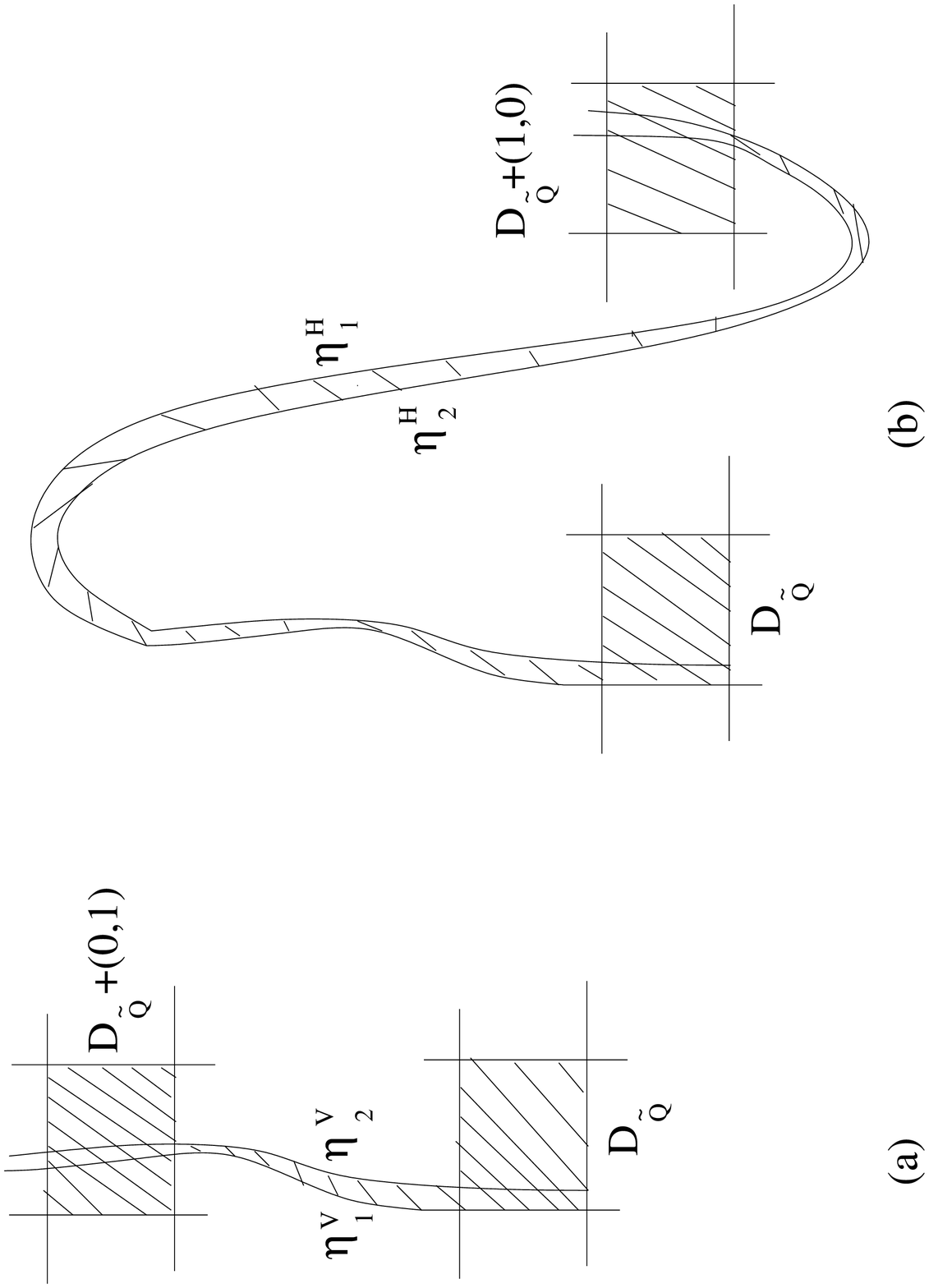}} 
\mbox{\includegraphics[width=13cm]{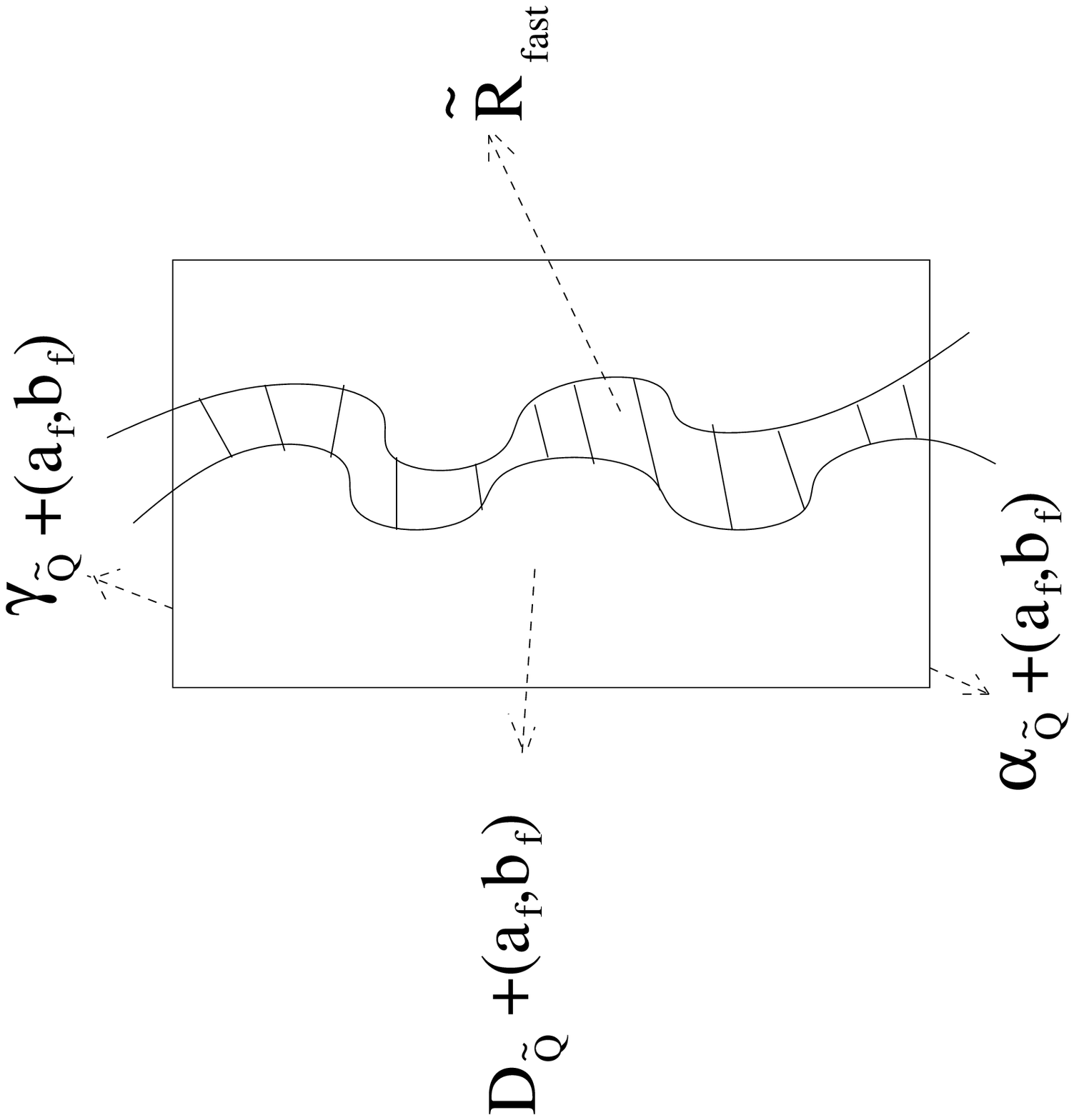}}
\end{center}

\end{document}